\documentclass[leqno]{article}

\begin{document}

\title{Ode to commutator operators}

\author{Stephen W.\ Semmes \\
	Rice University}

\date{}

\maketitle

\begin{abstract}
These brief remarks have been prepared in connection with a conference
in honor of my thesis advisor, Richard Rochberg.
\end{abstract}

\tableofcontents

\section{It was a dark and stormy night}
\label{reminiscence}
\setcounter{equation}{0}

	Or maybe it was just another pleasant springtime evening in
the suburbs of Stockholm.  The professor was sharing a house on the
institute grounds with a couple of messy bachelor mathematicians, and
numerous empty beer bottles were scattered about the room.  The
professor was enjoying a cool swedish brew too while his earnest
doctoral student listened intently.  ``People have been studying the
linear terms for fifty years'', or something like that.  This was in
the context of nonlinear problems in which functions are themselves
variables for some other functional.  Linearizations in simpler
situations have been studied for a much longer period of time.  Now it
would be time to pursue the quadratic and higher-order terms, as in
the calculus of Newton and Leibniz.

\section{Commutators}
\label{commutators}
\setcounter{equation}{0}

	If $\mathcal{A}$ is an associative algebra and $a, b \in
\mathcal{A}$, then the \emph{commutator} of $a, b \in \mathcal{A}$ is
defined by
\begin{equation}
	[a, b] = a \, b - b \, a.
\end{equation}
Commutators automatically satisfy the Leibniz rule for products,
\begin{eqnarray}
	[a, b \, c] = a \, b \, c - b \, c \, a 
	& = & a \, b \, c - b \, a \, c + b \, a \, c - b \, c \, a	\\
	& = & [a, b] \, c + b \, [a, c].			\nonumber
\end{eqnarray}
Thus a commutator is like a derivative.

	A typical situation of interest would be an algebra of linear
operators on a vector space of functions.  Suppose that $D$ is a
first-order differential operator acting on functions on some space,
which satisfies the ordinary Leibniz rule
\begin{equation}
	D(f_1, f_2) = (D f_1) \, f_2 + f_1 (D f_2)
\end{equation}
with respect to pointwise multiplication of functions.  Let $b$ be a
function on the same space, and let $M_b$ be the corresponding
operator of multiplication by $b$, $M_b(f) = b \, f$.  In this case,
the commutator of $D$ and $M_b$ is the same as the operator of
multiplication by $D b$.

\section{The Calder\'on commutator}
\label{calderon commutator}
\setcounter{equation}{0}

	Let $H$ be the Hilbert transform on the real line ${\bf R}$,
and let $D$ be the usual differentiation operator $d/dx$.  Suppose
that $A$ is a Lipschitz function on ${\bf R}$, so that the derivative
$a = A'$ of $A$ is bounded.  If $M_A$ again denotes the operator of
multiplication by $A$, then a celebrated theorem of Calder\'on
\cite{ca1} states that the commutator $[M_A, H D]$ determines a
bounded operator on $L^2({\bf R})$.  This would be trivial without the
Hilbert transform, since multiplication by a bounded function is a
bounded operator on $L^2$.

	Observe that
\begin{eqnarray}
	D \, [M_A, H] & = & D \, M_A \, H - D \, H \, M_A		\\
		& = & M_a \, H + M_A \, D \, H - D \, H \, M_A	\nonumber \\
		& = & M_a \, H + [M_A, H D].			\nonumber
\end{eqnarray}
Hence the boundedness of $[M_A, HD]$ and $D \, [M_A, H]$ on $L^2({\bf
R})$ are equivalent when $A$ is Lipschitz, since the Hilbert transform
is bounded too.  The boundedness of $D \, [M_A, H]$ on $L^2({\bf R})$
can be reformulated as the boundedness of $[M_A, H]$ as a mapping from
$L^2({\bf R})$ to the Sobolev space of functions on the real line with
first derivative in $L^2$.  This suggests many interesting questions
about analogous operators on other metric spaces, where it may be
easier to make sense of something like $|D f|$ than $D f$.

\section{Connes' noncommutative geometry}
\label{noncommutative geometry}
\setcounter{equation}{0}

	In Connes' theory \cite{cns}, the commutator $[T, M_b]$
between a singular integral operator $T$ and the operator $M_b$ of
multiplication by a function $b$ represents a sort of derivative of
$b$.  Situations in which $T^2 = I$ are of particular interest.  More
precisely, one may work with operators acting on vector-valued
functions here.  For instance, $T$ could be made up of Riesz
transforms.

	As analogues of classical integrals of products of functions,
Connes studies certain types of traces of products of commutators.
These products of commutators are normally not quite of trace class,
and so an extension of the trace due to Dixmier \cite{dxm} is
employed.  The one-dimensional case where such a commutator might be
of trace class is somewhat exceptional in this context.  The Dixmier
trace of a trace class operator is equal to $0$, and involves
asymptotic behavior of an operator in general.

	For that matter, ordinary derivatives involve limits.  This is
an important part of localization.  That is a natural feature to try
to have.

\section{Commutators as derivatives}
\label{commutators, derivatives}
\setcounter{equation}{0}

	Commutators often occur as derivatives of nonlinear
functionals.  This is a familiar theme in mathematics, where the
nonlinear functional involves some form of conjugation.  As a basic
example, the commutator $[M_b, T]$ of the operator $M_b$ of
multiplication by $b$ and another linear operator $T$ can be seen as
the derivative in $b$ at the origin of the conjugation $M_{e^b} \, T
\, M_{e^{-b}}$ of $T$ by multiplication by $e^b$.  This was studied in
\cite{c-r-w} in connection with weighted norm inequalities for
singular integral operators and bounded mean oscillation.

	The Calder\'on commutators correspond to derivatives of the
Cauchy integral operator on Lipschitz graphs, with respect to the
Lipschitz function being graphed.  As in \cite{c-m-2, c-m-3}, this is
closely related to the conjugation of the Hilbert transform by a
change of variables.  More precisely, the Cauchy integral operator can
be obtained as an analytic continuation of the conjugation of the
Hilbert transform by a change of variables, by interpreting it as a
conjugation of the Hilbert transform by a change of variables that
extends into the complex plane.

\section{Projections}
\label{projections}
\setcounter{equation}{0}

	Let $\mathcal{H}$ be a Hilbert space, and let $P$ be a bounded
linear operator on $\mathcal{H}$ which is a projection.  Thus $P^2 =
P$.  Let $V_0$ be the set of $v \in \mathcal{H}$ such that $P(v) = 0$,
and let $V_1$ be the set of $v \in \mathcal{H}$ such that $P(v) = v$.
These are closed linear subspaces of $\mathcal{H}$, and $P$ is the
projection of $\mathcal{H}$ onto $V_1$ with kernel $V_0$.  Consider $T
= 2 \, P - I$.  Equivalently, $T$ is characterized by the conditions
that $T(v) = -v$ when $v \in V_0$ and $T(v) = v$ when $v \in V_1$.  In
particular, $T^2 = I$.  Note that $T$ is self-adjoint if and only if
$P$ is, which happens exactly when $V_0$ and $V_1$ are orthogonal.  If
$P$ is the orthogonal projection of $L^2$ on the real line or unit
circle onto the Hardy space of functions of analytic type, then $T$
corresponds to the Hilbert transform.  Examples that are not
self-adjoint can be obtained from Cauchy integrals and using weights.

\section{Toeplitz operators}
\label{toeplitz}
\setcounter{equation}{0}

	Let $X$ be a space equipped with a topology and a Borel
measure, and perhaps a metric.  As a Hilbert space $\mathcal{H}$,
consider $L^2(X)$.  Suppose that $P$ is a bounded linear operator on
$\mathcal{H}$ which is a projection.  It might also be nice if $T = 2
\, P - I$ is something like a singular integral operator.  A key point
is for $[M_b, T]$ to be compact when $b$ is a continuous function on
$X$ which is continuous at infinity if $X$ is not compact.  Because
compact operators form a closed linear subspace of the space of
bounded linear operators on $\mathcal{H}$, it suffices to check this
for a dense class of functions $b$.  For instance, $[M_b, T]$ is more
regular when $b$ is a Lipschitz function with compact support.  Under
these conditions, one can consider Toeplitz operators $P \, M_b$ on
$V_1 = P(L^2(X))$.  These are the classical Toeplitz operators when
$X$ is the unit circle or real line and $P$ is the orthogonal
projection onto the associated Hardy space.

	In some cases, it may be convenient to allow vector-valued
functions on $X$.  For example, one might be interested in functions
on ${\bf R}^n$ with values in a Clifford algebra.  Using Clifford
analysis, one gets a Hardy space and a projection defined in terms of
singular integral operators.  In the classical situation, the product
of holomorphic functions is homolomorphic.  This does not work for
Clifford holomorphic functions.  Thus some of the usual structure is
not available.  One way to deal with this is to restrict one's
attention to scalar-valued functions $b$.  This ensures that $b$
commutes with elements of the Clifford algebra, so that the commutator
$[M_b, T]$ still behaves well.  Alternatively, one can be more careful
about the way in which a multiplication operator acts, from the left
or the right.  Because of noncommutativity, Clifford holomorphicity
can also be defined from the left or right.  The product of a Clifford
holomorphic function and a constant in the Clifford algebra is
Clifford holomorphic, when the constant is on the appropriate side.
Similarly, commutators of Clifford-valued multiplication operators and
singular integral operators also behave well when the operators
involve multiplication on opposite sides.

\section{Extensions}
\label{extensions}
\setcounter{equation}{0}

	Let $X$ be a compact metric space, and let $\mathcal{C}(X)$ be
the Banach algebra of continuous complex-valued functions on $X$ with
the supremum norm.  Let $\mathcal{H}$ be a complex Hilbert space, and
let $\mathcal{B}(\mathcal{H})$ be the Banach algebra of bounded linear
operators on $\mathcal{H}$.  Suppose that $\mathcal{E}$ is a
$C^*$-subalgebra of $\mathcal{B}(\mathcal{H})$ and that $\phi$ is a
$*$-algebra homomorphism of $\mathcal{E}$ onto $\mathcal{C}(X)$ whose
kernel is contained in the ideal $\mathcal{K}(\mathcal{H})$ of compact
operators on $\mathcal{H}$.  If $T \in \mathcal{E}$ and $\phi(T)$ is
an invertible continuous function on $X$, then $T$ is a Fredholm
operator on $\mathcal{H}$.  For if $R \in \mathcal{E}$ and $\phi(R) =
\phi(T)^{-1}$, then $R \, T - I$ and $T \, R - I$ are compact
operators.  It may be that $T$ is a compact perturbation of an
invertible operator.  However, this is not possible when the index of
$T$ is nonzero.

	As in \cite{b-d-f-1, b-d-f-2, b-d-f-3, d1}, these are
interesting circumstances to be in.  Classical Toeplitz operators in
one or more complex variables and pseudodifferential operators of
order $0$ provide basic examples of this type of situation.  From the
point of view of this theory, it is very natural to consider Toeplitz
operators associated to multiplication by scalar-valued functions in
the context of complexified Clifford analysis.  One can also compress
these Toeplitz operators to smaller spaces of Clifford-holomorphic
functions.

\section{Abstract elliptic operators}
\label{elliptic operators}
\setcounter{equation}{0}

	Let $X$ be a compact metric space equipped with a positive
Borel measure $\mu$.  As in \cite{a1}, it is interesting to look at
bounded linear operators on $L^2(X)$ whose commutators with
multiplication by continuous functions on $X$ are compact.  In
particular, Fredholm operators of this type generalize classical
elliptic pseudodifferential operators of order $0$ on compact smooth
manifolds.  One can also consider operators acting on vector-valued
functions, and other extensions, as explained in \cite{a1}.

	The Hilbert transform on the unit circle is a basic example.
Fractional integral operators of imaginary order are also examples.
These may be given as imaginary powers of unbounded nonnegative
self-adjoint operators, like Laplace operators.

\section{Fredholm indices}
\label{fredholm indices}
\setcounter{equation}{0}

	As in \cite{a1}, Fredholm operators with compact commutators
with operators of multiplication by continuous functions lead to other
Fredholm operators, using vector bundles on the underlying metric
space.  Thus one gets indices for all of these operators.  One also
gets Fredholm operators and hence indices associated to invertible
continuous complex-valued functions in the context of Section
\ref{extensions}.  This can be extended to matrix-valued functions,
since matrices of complex functions lead to matrices of operators
which can be interpreted as operators on other spaces.  Of course,
these two theories are closely related to each other.  In connection
with Toeplitz-type operators defined using Clifford analysis, one
should be careful about the precise domains and ranges of the
operators.  At any rate, nontrivial indices can be a good indication
of the presence of significant geometric or other structure.


\begin{thebibliography}{80}





\bibitem {a-r} J.~Anderson and R.~Rochberg, {\it Toeplitz operators
associated with subalgebras of the disk algebra}, Indiana University
Mathematics Journal {\bf 30} (1981), 813--820.

\bibitem {arv} W.~Arveson, {\it A Short Course on Spectral Theory},
Springer-Verlag, 2002.

\bibitem {a1} M.~Atiyah, {\it Global theory of elliptic operators}, in
{\it Proceedings of the International Conference on Functional
Analysis and Related Topics}, 21--30, University of Tokyo Press, 1970.

\bibitem {a2} M.~Atiyah, {\it $K$-Theory}, notes by D.~Anderson, 2nd
edition, Addison-Wesley, 1989.

\bibitem {bls} R.~Beals, {\it Characterization of pseudodifferential
operators and applications}, Duke Mathematical Journal {\bf 44}
(1977), 45--57.

\bibitem {b-d-s} F.~Brackx, R.~Delanghe, and F.~Sommen, {\it Clifford
Analysis}, Pitman, 1982.

\bibitem {b-d-f-1} L.~Brown, R.~Douglas, and P.~Fillmore, {\it
Extensions of $C^*$-algebras, operators with compact self-commutators,
and $K$-homology}, Bulletin of the American Mathematical Society {\bf
79} (1973), 973--978.

\bibitem {b-d-f-2} L.~Brown, R.~Douglas, and P.~Fillmore, {\it Unitary
equivalence modulo the compact operators and extensions of
$C^*$-algebras}, in {\it Proceedings of a Conference on Operator
Theory}, 58--128, Lecture Notes in Mathematics {\bf 345},
Springer-Verlag, 1973.

\bibitem {b-d-f-3} L.~Brown, R.~Douglas, and P.~Fillmore, {\it
Extensions of $C^*$-algebras and $K$-homology}, Annals of Mathematics
(2) {\bf 105} (1977), 265--324.

\bibitem {ca1} A.~Calder\'on, {\it Commutators of singular integral
operators}, Proceedings of the National Academy of Sciences U.S.A.\
{\bf 53} (1965), 1092--1099.

\bibitem {ca2} A.~Calder\'on, {\it Cauchy integrals on Lipschitz
curves and related operators}, Proceedings of the National Academy of
Sciences U.S.A.\ {\bf 74} (1977), 1324--1327.

\bibitem {ca3} A.~Calder\'on, {\it Commutators, singular integrals on
Lipschitz curves and applications}, in {\it Proceedings of the
International Congress of Mathematicians (Helsinki, 1978)}, 85--96,
Academiae Scientiarum Fennicae, 1980.

\bibitem {c-m-m} R.~Coifman, A.~McIntosh, and Y.~Meyer, {\it
L'Integrale de Cauchy d\'efinit un op\'erateur born\'e sur $L^2$ pour
les courbes lipschitziennes}, Annals of Mathematics (2) {\bf 116}
(1982), 361--387.

\bibitem {c-m-1} R.~Coifman and Y.~Meyer, {\it Au-del\`a des
Op\'erateurs Pseudo-Diff\'erentiels}, Ast\'erisque {\bf 57}, 1978.

\bibitem {c-m-2} R.~Coifman and Y.~Meyer, {\it Le th\'eor\`eme de
Calder\'on par les ``m\'ethodes de variables r\'eelle''}, Comptes
Rendus Hebdomadaires des S\'eances de l'Acad\'emie des Sciences Paris
S\'eries A et B {\bf 289} (1979), A425--A428.

\bibitem {c-m-3} R.~Coifman and Y.~Meyer, {\it Une g\'en\'eralisation
du th\'eor\`eme de Calder\'on sur l'int\'egrale de Cauchy}, in {\it
Fourier Analysis}, 87--116, Asociacion Matematica Espa\~{n}ola, 1980.

\bibitem {c-r} R.~Coifman and R.~Rochberg, {\it Projections in
weighted spaces, skew projections and inversion of Toeplitz
operators}, Integral Equations and Operator Theory {\bf 5} (1982),
145--159.

\bibitem {c-r-w} R.~Coifman, R.~Rochberg, and G.~Weiss, {\it
Factorization theorems for Hardy spaces in several variables}, Annals
of Mathematics (2) {\bf 103} (1976), 611--635.

\bibitem {c-w-1} R.~Coifman and G.~Weiss, {\it Analyse Harmonique
Non-Commutative sur Certains Espaces Homog\`enes: \'Etude de Certained
Int\'egrales Singuli\`eres}, Lecture Notes in Mathematics {\bf 242},
1971.

\bibitem {c-w-2} R.~Coifman and G.~Weiss, {\it Extensions of Hardy
spaces and their use in analysis}, Bulletin of the American
Mathematical Society {\bf 83} (1977), 569--645.

\bibitem {cns} A.~Connes, {\it Noncommutative Geometry}, Academic
Press, 1994.

\bibitem {c-s-t-1} A.~Connes, D.~Sullivan, and N.~Teleman, {\it
Formules locales pour les classes de Pontrjagin topologiques}, Comptes
Rendus de l'Acad\'emie des Sciences Paris S\'erie I Mathematique {\bf
317} (1993), 521--526.

\bibitem {c-s-t-2} A.~Connes, D.~Sullivan, and N.~Teleman, {\it
Quasiconformal mappings, operators on Hilbert space, and local
formulae for characteristic classes}, Topology {\bf 33} (1994),
663--681.

\bibitem {c} M.~Cotlar, {\it Convolution Operators and Factorization},
McGill University, 1972.

\bibitem {c-j-m-r} M.~Cwikel, B.~Jawerth, M.~Milman, and R.~Rochberg,
{\it Differential estimates and commutators in interpolation theory},
in {\it Analysis at Urbana}, Volume II, 170--220, Cambridge University
Press, 1989.

\bibitem {c-k-m-r} M.~Cwikel, N.~Kalton, M.~Milman, and R.~Rochberg,
{\it A unified theory of commutator estimates for a class of
interpolation methods}, Advances in Mathematics {\bf 169} (2002),
241--312.

\bibitem {dxm} J.~Dixmier, {\it Existence de traces non normales},
Comptes Rendus Hebdomadaires des Sc\'eances de l'Acad\'emie des
Sciences Paris S\'eries A et B {\bf 262} (1966), A1107--A1108.

\bibitem {d1} R.~Douglas, {\it $C^*$-Algebra Extensions and
$K$-Homology}, Princeton University Press, 1980.

\bibitem {d2} R.~Douglas, {\it Banach Algebra Techniques in Operator
Theory}, 2nd edition, Springer-Verlag, 1998.

\bibitem {drn} P.~Duren, {\it Theory of $H^p$ Spaces}, Academic Press,
1970.

\bibitem {drn-s} P.~Duren and A.~Schuster, {\it Bergman Spaces},
American Mathematical Society, 2004.

\bibitem {f-r} M.~Feldman and R.~Rochberg, {\it Singular value
estimates for commutators and Hankel operators on the unit ball and
the Heisenberg group}, in {\it Analysis and Partial Differential
Equations}, 121--159, Dekker, 1990.

\bibitem {g} J.~Garnett, {\it Bounded Analytic Functions}, revised
edition, Springer-Verlag, 2007.

\bibitem {g-m} J.~Gilbert and M.~Murray, {\it Clifford Algebras and
Dirac Operators in Harmonic Analysis}, Cambridge University Press,
1991.

\bibitem {h-k-z} H.~Hedenmalm, B.~Korenblum, and K.~Zhu, {\it Theory
of Bergman Spaces}, Springer-Verlag, 2000.

\bibitem {hfm} K.~Hoffman, {\it Banach Spaces of Analytic Functions},
Dover, 1988.

\bibitem {j} S.~Janson, {\it Mean oscillation and commutators of
singular integral operators}, Arkiv f\"or Matematik {\bf 16} (1978),
263--270.

\bibitem {j-w} S.~Janson and T.~Wolff, {\it Schatten classes and
commutators of singular integral operators}, Arkiv f\"or Matematik
{\bf 20} (1982), 301--310.

\bibitem {j-r-w} B.~Jawerth, R.~Rochberg, and G.~Weiss, {\it
Commutator and other second order estimates in real interpolation
theory}, Arkiv f\"or Matematik {\bf 24} (1986), 191--219.

\bibitem {jl} J.-L.~Journ\'e, {\it Calder\'on--Zygmund Operators,
Pseudodifferential Operators, and the Cauchy Integral of Calder\'on},
Lecture Notes in Mathematics {\bf 994}, Springer-Verlag, 1983.

\bibitem {k} J.~Kigami, {\it Analysis on Fractals}, Cambridge
University Press, 2001.

\bibitem {kss} P.~Koosis, {\it Introduction to $H_p$ Spaces}, 2nd
edition, with two appendices by V.~Havin, Cambridge University Press,
1998.

\bibitem {sk} S.~Krantz, {\it Function Theory of Several Complex
Variables}, AMS Chelsea, 2001.

\bibitem {k-l-1} S.~Krantz and S.-Y.~Li, {\it Boundedness and
compactness of integral operators on spaces of homogeneous type and
applications, I}, Journal of Mathematical Analysis and Applications
{\bf 258} (2001), 629--641.

\bibitem {k-l-2} S.~Krantz and S.-Y.~Li, {\it Boundedness and
compactness of integral operators on spaces of homogeneous type and
applications, II}, Journal of Mathematical Analysis and Applications
{\bf 258} (2001), 642--657.

\bibitem {l-m-q} C.~Li, A.~McIntosh, and T.~Quin, {\it Clifford
algebras, Fourier transforms, and singular convolution operators on
Lipschitz surfaces}, Revista Matem\'atica Iberoamericana {\bf 10}
(1994), 665--721.

\bibitem {l-m-s} C.~Li, A.~McIntosh, and S.~Semmes, {\it Convolution
singular integrals on Lipschitz surfaces}, Journal of the American
Mathematical Society {\bf 5} (1992), 455--481.

\bibitem {m-m-v} P.~Mattila, M.~Melnikov, and J.~Verdera, {\it The
Cauchy integral, analytic capacity, and uniform rectifiability},
Annals of Mathematics (2) {\bf 144} (1996), 127--136.

\bibitem {mel} M.~Melnikov, {\it Analytic capacity: A discrete
approach and the curvature of measure}, Mathematics Sbornik {\bf 186}
(1995), 827--846.

\bibitem {mel-v} M.~Melnikov and J.~Verdera, {\it A geometric proof of
the $L^2$ boundedness of the Cauchy integral on Lipschitz curves},
International Mathematics Research Notices {\bf 1995}, 325--331.

\bibitem {na1} A.~Nahmod, {\it Geometry of operators and spectral
analysis on spaces of homogeneous type}, Comptes Rendus de
l'Acad\'emie des Sciences Paris S\'erie I Math\'ematique {\bf 313}
(1991), 721--725.

\bibitem {na2} A.~Nahmod, {\it Generalized uncertainty principles on
spaces of homogeneous type}, Journal of Functional Analysis {\bf 119}
(1994), 171--209.

\bibitem {n1} N.~Nikolski, {\it Treatise on the Shift Operator}, with
an appendix by S.~Hru{\v s}{\v c}ev and V.~Peller, translated from the
Russian by J.~Peetre, Springer-Verlag, 1986.

\bibitem {n2} N.~Nikolski, {\it Operators, Functions, and Systems: An
Easy Reading}, Volume 1, {\it Hardy, Hankel, and Toeplitz}, Volume 2,
{\it Model Operators and Systems}, translated from the French by
A.~Hartmann and revised by the author, American Mathematical Society,
2002.

\bibitem {pt} J.~Partington, {\it An Introduction to Hankel
Operators}, Cambridge University Press, 1988.

\bibitem {p-r} J.~Peetre and R.~Rochberg, {\it Higher order Hankel
forms}, in {\it Multivariable Operator Theory}, 283--306, American
Mathematical Society, 1995.

\bibitem {p} V.~Peller, {\it Hankel Operators and their Applications},
Springer-Verlag, 2003.

\bibitem {pr} S.~Power, {\it Hankel Operators on Hilbert Space},
Pitman, 1982.

\bibitem {r1} R.~Rochberg, {\it Toeplitz operators on weighted $H^p$
spaces}, Indiana University Mathematics Journal {\bf 26} (1977),
291--298.

\bibitem {r2} R.~Rochberg, {\it Trace ideal criteria for Hankel
operators and commutators}, Indiana University Mathematics Journal
{\bf 31} (1982), 913--925.

\bibitem {r3} R.~Rochberg, {\it Higher order estimates in complex
interpolation theory}, Pacific Journal of Mathematics {\bf 174}
(1996), 247--267.

\bibitem {r4} R.~Rochberg, {\it Higher-order Hankel forms and
commutators}, in {\it Holomorphic Spaces}, 155--178, Cambridge
University Press, 1998.

\bibitem {r-s-1} R.~Rochberg and S.~Semmes, {\it A decomposition
theorem for BMO and applications}, Journal of Functional Analysis {\bf
67} (1986), 228--263.

\bibitem {r-s-2} R.~Rochberg and S.~Semmes, {\it End point results for
estimates of singular values of singular integral operators}, in {\it
Contributions to Operator Theory and its Applications}, 217--231,
Birkh\"auser, 1988.

\bibitem {r-s-3} R.~Rochberg and S.~Semmes, {\it Nearly weakly
orthonormal sequences, singular value estimates, and Calder\'on--Zygmund
operators}, Journal of Functional Analysis {\bf 86} (1989), 237--306.

\bibitem {r-w} R.~Rochberg and G.~Weiss, {\it Derivatives of analytic
families of Banach spaces}, Annals of Mathematics (2) {\bf 118}
(1983), 315--347.

\bibitem {r} W.~Rudin, {\it Function Theory in the Unit Ball of ${\bf
C}^n$}, Springer-Verlag, 1980.

\bibitem {sar} D.~Sarason, {\it Function Theory on the Unit Circle},
Department of Mathematics, Virginia Polytechnic Institute and State
University, 1978.

\bibitem {sng} I.~Singer, {\it Future extensions of index theory and
elliptic operators}, in {\it Prospects in Mathematics}, 171--185,
Princeton University Press, 1971.

\bibitem {st1} E.~Stein, {\it Singular Integrals and Differentiability
Properties of Functions}, Princeton University Press, 1970.

\bibitem {st2} E.~Stein, {\it Harmonic Analysis: Real-Variable
Methods, Orthogonality, and Oscillatory Integrals}, with the
assistance of T.~Murphy, Princeton University Press, 1993.

\bibitem {s-w} E.~Stein and G.~Weiss, {\it Introduction to Fourier
Analysis on Euclidean Spaces}, Princeton University Press, 1971.

\bibitem {sz1} R.~Strichartz, {\it Analysis on fractals}, Notices of
the American Mathematical Society {\bf 46} (1999), 1199--1208.

\bibitem {sz2} R.~Strichartz, {\it Differential Equations on
Fractals: A Tutorial}, Princeton University Press, 2006.

\bibitem {u} A.~Uchiyama, {\it On the compactness of operators of
Hankel type}, T\^ohoku Mathematical Journal (2) {\bf 30} (1978),
163--171.

\bibitem {v} J.~Verdera, {\it $L^2$ Boundedness of the Cauchy integral
and Menger curvature}, in {\it Harmonic Analysis and Boundary Value
Problems}, 139--158, American Mathematical Society, 2001.

\bibitem {w} N.~Weaver, {\it Lipschitz Algebras}, World Scientific,
1999.

\bibitem {z1} K.~Zhu, {\it Operator Theory in Function Spaces},
Dekker, 1990.

\bibitem {z2} K.~Zhu, {\it Spaces of Holomorphic Functions in the Unit
Ball}, Springer-Verlag, 2005.


\end{thebibliography}
\end{document}